\documentclass{article}

\usepackage{amssymb}
\usepackage{amsmath}

\let\<\langle
\let\>\rangle
\font\bigcmsy=cmsy10.pk scaled 2000
\def\bigtimes{\mathop{\,\vrule width0pt depth2pt height8pt
            \smash{\lower2pt\hbox{\bigcmsy\char'002}}\,}\limits}

\begin{document}


\begin{center}
\Large{\textbf{Multielement order separability in free products of
groups}} \footnote{This work was supported by the Russian Foundation
for Basic Research, project no.~05-01-00895}
\end{center}

\begin{center}
\textbf{Vladimir V. Yedynak }
\end{center}

\begin{abstract}
We introduce the notion of multielement order separability and study this property for free groups and free
products.


\end{abstract}

\section{Introduction.}

Definition 1. We say that a group $G$ is \textsl{$n$-element order
separable} if, for each elements $u_1, ...,\ u_n$ of $G$ such that
$u_i$ and $u_j^{\pm1}$ are not conjugate for $i\ne j$, there exists
a homomorphism of $G$ onto a finite group such that the orders of
the images of $u_1, ..., u_n$ are pairwise different.

It is easy to see that, if a free product of two groups is
$n$-element order separable, then the free factors are $k$-element
order separable for each $k\leqslant n$.

The converse is not true for $n>2$. For example, the infinite
dihedral group $\<a\>_2*\<b\>_2$ is not 3-element order separable.

So, we require additionally each two nontrivial elements from
different free factors not to have equal finite orders. It is known
that free groups are 2-element order separable (Equations over
groups, quasivarieties, and a residual property of a free group,
\textsl{Journal of group theory} \textbf{2}: 319-327), and,
moreover, this property is inherited by free products (Separability
with respect to order, \textsl{Vestnik Mosk. Univ. Ser. I Mat.
Mekh.} \textbf{3}: 56-58).

It was proven that free metabelian groups are 2-element order
separable, but this property fails for free 3-step solvable groups
(Elements with the same normal closure in a metabelian group,
\textsl{The Quarterly Journal of Mathematics} \textbf{58(1)}:
23-29).

\textbf{ Theorem 1.}\textsl{ Any absolutely free group is
$n$-element order separable for each $n$.}

\textbf{ Theorem 2.}\textsl{ Suppose that groups $A$ and $B$ are
periodic and $n$-element order separable for each $n$. If $|a|\neq
|b|$ for any $a\in A\setminus\{ e\}, b\in B\setminus\{ e\}$, then
the group $A\ast B$ is also $n$-element order separable for each
$n$.}

\textbf{ Theorem 3.}\textsl{ If groups $A$ and $B$ are $n$-element
order separable for each $n$ and $|a|\neq |b|$ for any $a\in
A\setminus\{ e\}, b\in B\setminus\{ e\}$, such that $|a|<\infty,
|b|<\infty$, then the free product $G=A\ast B$ is $3$-element order
separable.}

\section{Notations and definitions.}

Suppose that $G=A\ast B$, where $A$ and $B$ are finite groups.
Consider an oriented graph $\Gamma$ with the following properties:

(1) For each vertex $p$ and each $c\in A\cup B\setminus\{1\}$, there
are exactly one oriented edge labelled by $c$ and starting at $p$,
exactly one oriented edge labelled by $c$ and ending at $p$.

(2) For each vertex $p$, the maximal oriented connected subgraph
$A(p)$ containing $p$, whose edges are labelled by elements of the
group $A$, is the Cayley graph of the group $A$ with respect to the
set of generators $A$. Similarly, for each vertex $p$, the maximal
oriented connected subgraph $B(p)$ containing $p$, whose edges are
labelled by elements from the group $B$, is the Cayley graph of the
group $B$ with respect to the set of generators $B$.

If $e$ is an edge of the graph $\Gamma$, then Lab($e$) denotes the
label of $e$. If $s=e_1...e_n$ is a path in $\Gamma$, then Lab($s$)
= Lab($e_1)...$Lab($e_n$) denotes the label of the path $s$,
$\alpha(s)=\alpha(e_1)$ is the beginning of $e_1$,
$\omega(s)=\omega(e_n)$ is the end of $e_n$.

Definition 2. \textsl{Let $x\in G$ be a cyclically reduced element
not belonging to $A\cup B$. We say that a closed path $e_1...e_n$
whose label equals $x^k$ is an $x$-cycle if
Lab($e_{il+1}...e_{il+l})$ is the normal form of the element $x$,
where $l$ is the length of the word $x$, $i$ is an arbitrary
positive integer, and indices are modulo $n$. The number $k$ is
called the length of this $x$-cycle.}

Definition 3. \textsl{If two different edges of a cycle lies in the
same subgraph $A(p)$ or $B(p)$, we say that this cycle has close
edges.}

Let $\Gamma$ be a graph with Properties (1), (2). Take a set of
vertices $S$ in the graph $\Gamma$. Consider an arbitrary function
$f: S\rightarrow\{A, B\}$. Consider $t$ copies $\Delta_1, ...,
\Delta_t$ of $\Gamma$. Let $v_i\in\Delta_i$ be the vertex
corresponding to a vertex $v\in\Gamma$, let $q_i\in\Delta_i$ be the
edge corresponding to an edge $q\in\Gamma$ and let $S_i$ be the set
of vertices of $\Delta_i$ corresponding to the set $S$. For each
$i\in\{1,\dots,t\}$, for each edge $q$ starting or ending at a
vertex from $S$ and each $p\in S$, we delete edges labelled by the
elements of the group $f(p)$ whose end points coincide with $p_i$,
if $q_i$ starts at $p_i$, ends at $v_i$ and Lab$(q_i)\in f(p)$
connect the vertices $p_i$ and $v_{i+1}$ (subscripts modulo $t$) by
an oriented edge starting at $p_i$ whose labell equals Lab$(q_i)$,
if $q_i$ starts at $v_i$, ends at $p_i$ and Lab$(q_i)\in f(p)$
connect the vertices $p_i$ and $v_{i+1}$ (subscripts modulo $t$) by
an oriented edge starting at $v_{i+1}$ whose labell equals
Lab$(q_i)$. Let $S=\{s_1,\dots, s_k\}$. The obtained graph
$\gamma_t(\Gamma; s_1,\dots, s_k; f(s_1),\dots, f(s_k))$ satisfies
Properties (1), (2).

The group $G$ acts on the right on the set of vertices of the graph
$\Gamma$ with Properties (1), (2) by the following way. If $a\in A$
and $p$ is a vertex of $\Gamma$ then there exists exactly one edge
$h$ starting at $p$ and labelled by $a$. Put $p\cdot a=\omega(h)$.
(Here and in what follows, $\alpha(h)$ and $\omega(h)$ denote the
start and end points of an oriented edge $h$.) The action of group
$B$ is defined similarly.

In what follows, having a homomorphism $\varphi$ of the group $A\ast
B$ we always assume that the Cayley graph of the group
$\varphi(A\ast B)$ is constructed with respect to the set of
generators $\varphi(A\cup B)$.

\section{Auxiliary Lemmas.}

\textbf{ Lemma 1.}\textsl{ If $G=A\ast B$, where the groups $A$ and
$B$ are finite, $w_1, ..., w_k\in C\setminus e$ (throughout this
paper, $C$ denotes the Cartesian subgroup of the free product of
groups $A$ and $B$), and $p$ is a prime number, $N$ is a positive
integer, then there exists a homomorphism $\varphi$ of $G$ onto a
finite group such that, for each $i$, $\varphi(w_i)$ is a
$p$-element and $|\varphi(w_i)|>p^N$.}

\textsl{ Proof.}

Since $C$ is a free group, there exists a homomorphism $\varphi_1$
of $C$ onto a finite $p$-group such that $\varphi_1(w_i^j)\neq 1$
for $i=1, ..., k, j=1, ..., N$ (Foundations of group theory,
\textsl{Nauka}, theorem 14.2.2). Put $N=$ Ker $\varphi_1$. Consider
the subgroup $N'=\cap_{g\in G}(g^{-1}Ng)$. Then $N'\vartriangleleft
G, |G : N'|< \infty$ and the order of the image of $w_i$ in the
group $G/N'$ is a power of $p$ which is greater than $N$ for each
$i$. The natural homomorphism $\varphi : G\longrightarrow G/N'$ is
as required.

\textbf{ Lemma 2.} \textsl{ Suppose that $G=A\ast B$, where $A$ and
$B$ are finite groups. Let $S$ be the finite set of elements, such
that for any $u\in S$ the following conditions are satisfied: $u$ is
a nonunit cyclically reduced element of the Cartesian subgroup $C$,
furthermore $u=u'^m$, where $u'$ is not a proper power and if $m>1$
then for each $m'\neq m$ such that $m'|m\ u'^{m'}\notin C$. Then for
any prime $p$ there exists a homomorphism $\varphi$ of $G$ onto a
finite group such that for each $u\in S$ $\varphi(u)$ is a nonunit
$p$-element, and all $u$-cycles in the Cayley graph of $\varphi(G)$
have no close edges.}

\textsl{ Proof.}

Consider a prime $p$. Fix an element $u$ of $S$. Put $u'=y_1...y_n$
--- the normal form of $u'$. Fix $\mu, \nu$ such that $1\leqslant
\mu, \nu\leqslant n$. Fix also $z\in A\cup B$. We consider that if
$z=1$ then $\mu\neq\nu+1$ and $\nu-\mu<n-1$.

If $\mu>1$ then $\mu'=y_{\mu}...y_n$ otherwise $\mu'=1$. If $\nu<n$
then $\nu'=y_1...y_{\nu}$ otherwise $\nu'=1$. Consider an element
$\chi_{z, \mu, \nu}=\nu'z\mu'$. It is easy to see that $\chi_{z,
\mu, \nu}\notin\<u'\>$. Let us to prove that there exists a
homomorphism $\varphi_{z, \mu, \nu}$ of $G$ onto a finite group such
that for each $w\in S$ $\varphi_{z, \mu, \nu}(w)$ --- nonunit
$p$-element, $\varphi_{z, \mu, \nu}(\chi_{z, \mu,
\nu})\notin\<\varphi_{z, \mu, \nu}(u')\>$, and $|\varphi_{z, \mu,
\nu}(u')|=m|\varphi_{z, \mu, \nu}(u)|$. Lemma 1 shows that there
exists a homomorphism $\varphi_1$ of $G$ onto a finite group such
that $\varphi_1(w)$ is a nonunit $p$-element for each $w\in S$. It
is obvious that $|\varphi_1(u')|=m'|\varphi_1(u)|$, where $m'$ is a
divisor of $m$. Since $u'^m\in C$ and $u'^{m'|\varphi_1(u)|}\in$ Ker
$\varphi_1\subset C$ then $u'^{m'}\in C$. Hence we deduce that
$m'=m$ and $|\varphi_1(u')|=m|\varphi_1(u)|$. We may assume that
$\varphi_1(\chi_{z, \mu, \nu})\in\<\varphi_1(u')\>$.

Let $\Gamma$ be the Cayley graph of the group $\varphi_1(G)$.
Consider an arbitrary vertex $r$ in the graph $\Gamma$ and a path
$L=f_1...f_t$ starting at $r$ such that Lab($L$) is the normal form
of $\chi_{z, \mu, \nu}$. Suppose that $\omega(L)=q$. By the
assumption, there exists a $u'$-cycle $S=e_1...e_k$ starting at $q$,
and there exists an edge $e_l$ such that $\omega(e_l)=\alpha(L)$ and
Lab($e_1...e_l)=u'^v$. Let $\Delta_1, ..., \Delta_p$ be copies of
the graph $\Gamma$ and let $e_i^j$ be the edge of $\Delta_j$
corresponding to the edge $e_i$ in $\Gamma$. Let $\chi_{z, \mu,
\nu}=x_1...x_s,$ be the normal form of $\chi_{z, \mu, \nu}$.

If $y_n\in A$ then $D$ denotes $A$ and $F$ denotes $B$, if $y_n\in
B$ then $D$ denotes $B$, and $F$ denotes $A$.

If $x_1$ and $y_n$ belong to different free factors, we construct
the graph $\Delta=\gamma_p(\Gamma; r, \omega(e_{l+1}); D, D)$ from
the copies $\Delta_i$.

If $x_1$ and $y_n$ belong to the same free factor, then we construct
the graph $\Delta=\gamma_p(\Gamma; r, \alpha(e_l); F,F)$ from the
copies $\Delta_i$.

It is obvious that the length of each $u'$-cycle equals either
$|\varphi_1(u')|$ or $p|\varphi_1(u')|$ in $\Delta$. Besides, the
length of each $w$-cycle equals either $|\varphi_1(w)|$ or
$p|\varphi_1(w)|$ in $\Delta$ for any $w\in S$. The group $G$ acts
on the right on the set of vertices of $\Delta$. Besides, if $r_i$
is a vertex of $\Delta_i$ corresponding to the vertex $r$ of
$\Gamma$, then $r_i\cdot \chi_{z, \mu, \nu}u^l \neq r_i$ in $\Delta$
for any integer $l$. Thus, we found the required homomorphism
$\varphi_{z, \mu, \nu}$. Besides, $|\varphi_{z, \mu,
\nu}(u')|=m|\varphi_{z, \mu, \nu}(u)|$ and $\varphi_{z, \mu,
\nu}(w)$ --- nonunit $p$-element for any $w\in S$.

Consider a homomorphism $\varphi_u: G\rightarrow\prod_{z, \mu,
\nu}(\varphi_{z, \mu, \nu}(G))$ defined by the formula:
$g\mapsto\prod_{z, \mu, \nu}\varphi_{z, \mu, \nu}(g)$, where $z\in
A\cup B, 1\leqslant\mu, \nu\leqslant n$ and if $z=1$ then
$\mu\neq\nu+1, \nu-\mu<n-1$. In the Cayley graph $Cay(\varphi_u(G))$
of the group $\varphi_u(G)$ all $u'$-cycles have no close edges and,
$\varphi_u(w)$ --- is a nonunit $p$-element for any $w\in S$.
Furthermore $|\varphi_u(u')|$ is the less common multiple of numbers
$m|\varphi_{z, \mu, \nu}(u)|$ that is $m|\varphi_u(u)|$. Hence all
$u$-cycles in $Cay(\varphi_u(G))$ also have no close edges. It is
obvious that the homomorphism $\varphi: G\rightarrow\prod_{u\in
S}(\varphi_u(G))$ defined by the formula $g\mapsto\prod_{u\in
S}\varphi_u(g)$ is as required.

\textbf{ Lemma 3.}\textsl{ Suppose that $G=A\ast B$, where $A$ and
$B$ are finite groups, $u_1, ..., u_m$ are cyclically reduced
elements of the Cartesian subgroup such that each two of them do not
belong to conjugate cyclic subgroups, $\pi$ is a finite set of prime
numbers. Besides $u_i=u_i'^{m_i}$, where $u_i'$ is not a proper
power and if $m_i>1$ then for each $m_i'\neq m_i$ such that
$m_i'|m_i\ u_i'^{m_i'}\notin C$. Then there exists a homomorphism
$\varphi$ of $G$ onto a finite group such that the images of $u_1,
..., u_m$ have pairwise different orders and $p \nmid
|\varphi(u_i)|$ for each $i$ and for each $p\in \pi$.}

\textsl{ Proof. }

Let us prove by the induction on $m$. If $m=1$, the assertion
follows from Lemma 1.

Let us prove the following statement. For elements $u_1, ..., u_m$,
there exists a homomorphism $\varphi_1$ of $G$ onto a finite group
such that the images of $u_1, ..., u_m$ are nonidentity
$p$-elements, $p\notin\pi$; in addition $\{u_1,...,
u_m\}=\alpha\cup\beta$, where $\alpha$ and $\beta$ are nonempty sets
with empty intersection, moreover, for each $u_i\in\alpha$ and
$u_j\in\beta$, we have $|\varphi_1(u_i)|
> |\varphi_1(u_j)|$.

According to Lemma 2, there exists a homomorphism $\psi$ of $G$ onto
a finite group such that $\psi(u_1),...,\psi(u_m)$ are nonidentity
$p$-elements, where $p$ is a prime and $p\notin\pi$. Furthermore, if
$\Gamma$ is the Cayley graph of the group $\psi(G)$, then, in the
graph $\Gamma$, all $u_i$-cycles have no close edges for each $i$.
We may assume that $|\psi(u_1)|= ... = |\psi(u_m)|$. Consider some
set of nonnegative integers $k$. For each $k$ from this set, we
construct a graph $\Gamma_k$ satisfying (1), (2), and the following
properties:

(3) the lengths of maximal $u_i$-cycles coincide for $i=1,\dots,m$;

(4) for each $i$, the length of each $u_i$-cycle divides the length
of a maximal $u_i$-cycle;

(5) $u_i$-cycles have no close edges;

(6) if $k>0$, then there exists a path $R_k$ of length $k$ lying in
some maximal $u_1$-cycle and in all maximal $u_j$-cycles for each
$j=2,...,m$;

Properties (1), (2) and (3)--(6) are obviously true for the graph
$\Gamma_0=\Gamma$.

Having a graph $\Gamma_k$, we construct a new graph
$\Gamma^{\ast}_{k+1}$ satisfying Properties (1), (2) and (4), (5).
Hence we shall deduce that this graph satisfies Properties (3), (6)
so we put $\Gamma_{k+1}=\Gamma^{\ast}_{k+1}$; otherwise, owing to
the Property (5), we will prove that the least common multiple of
the lengths of all $u_1$-cycles and the least common multiple of the
lengths of all $u_j$-cycles in $\Gamma_{k+1}^{\ast}$ do not coincide
for some $j$.

Let $n$ be the length of a maximal $u_1$-cycle in $\Gamma_k$.
Consider a $u_1$-cycle $S$ in $\Gamma_k$ of length $n$. Suppose that
if $k>0$ then $S$ contains $R_k$.

Suppose that $k=0$. Let $S=f_1...f_k$. Put $s=\omega(f_1)$ and
$t=\omega(f_2)$. If Lab$(f_1)\in A$ then $D$ denotes $A$ and $F$
denotes $B$, if Lab$(f_1)\in B$ then $D$ denotes $B$ and $F$ denotes
$A$. Consider the graph
$\Gamma_1^{\ast}=\gamma_{n^2}(\gamma_n(\Gamma_0; s; D); t_{2,1};
F)$. Let $t_2\in \gamma_n(\Gamma_0; s; D)$ be the vertex of the
second copy of $\Gamma_0$ corresponding to $t$. Then
$t_{2,1}\in\Gamma_1^{\ast}$ is the vertex corresponding to $t_2$,
which belongs the first copy of the $\gamma_n(\Gamma_0; s; D)$. Put
$R_1'=f_2'$, where $f_2'$ is the edge of the graph $\Gamma_1^{\ast}$
labelled Lab($f_2$) and starting at the vertex $t_{2,1}$.

Consider the case $k>0$.  Suppose that $f$ is the edge of $S$ next
to $\omega(R_k)$, where $S$ is a maximal $u_1$-cycle containing
$R_k$. If Lab$(f)\in A$ then $D$ denotes $A$, if Lab$(f)\in B$ then
$D$ denotes $B$. Put $q=\omega(f)$. Consider the graph
$\Gamma^{\ast}_{k+1}=\gamma_n(\Gamma_k; q; D)$. Put
$R_{k+1}'=R_{k,1}\cup f_1$, where the path $R_{k,1}$ and the edge
$f_1$ belong to the first copy of $\Gamma_k$ in
$\Gamma^{\ast}_{k+1}$ and correspond to $R_k$ and $f$.

Consider a graph $X$, satisfying Properties (1),(2). Suppose that
$p$ is a vertex of $X$ and each $u_1$-cycle in $X$ has no close
edges, then, in the graph $\gamma_r(X; p; K)$, where $K$ is $A$ or
$B$, each $u_1$-cycle $V$ is the union of $r$ copies of a
$u_1$-cycle $W$ of the graph $X$. This means that, if $W$ has no
close edges, then $V$ has no close edges too. Let l($V$) and l($W$)
be the lengths of the $u_1$-cycles $V$ and $W$ respectively. Since
$V$ and $W$ have no close edges l($V$)=$r$l($W$). Owing to these
properties, we conclude that the graph $\Gamma^{\ast}_{k+1}$
satisfies Property (5), and since $\Gamma_{k+1}^{\ast}=\gamma_n(Y;z;
K)$, $n$ and length of each $u_i$-cycle of the graph $Y$ for
$i=1,\dots,m$ are powers of $p$ Property (4) is also fulfilled for
the graph $\Gamma_{k+1}^{\ast}$.

The length of $S$ becomes $n$ times greater, and there exists a
$u_j$-cycles in $\Gamma_k$, whose length also becomes $n$ times
greater in $\Gamma^{\ast}_{k+1}$ for each $j=2,...,m$. Properties
(4) and (5) shows that, if there exists a $u_j$-cycle in
$\Gamma^{\ast}_{k+1}$, whose length coincide with the length of a
$u_1$-cycle which is maximal in $\Gamma^{\ast}_{k+1}$, then this
$u_j$-cycle must contain the path $R_{k+1}'$. So if we suppose that
Condition (3) is fulfilled for $\Gamma^{\ast}_{k+1}$, then
$R'_{k+1}$ is contained in all maximal $u_j$-cycles of
$\Gamma^{\ast}_{k+1}$ and, if this is true for each $j=2,...,m$,
then this graph satisfies Property (6), and we put
$\Gamma_{k+1}=\Gamma^{\ast}_{k+1}$ and $R_{k+1}=R'_{k+1}$. Since
$u_1, u_j$ do not belong to conjugate cyclic subgroups for each
$j=2,...,m$, there exists a number $k$ such that the graph
$\Gamma^{\ast}_{k+1}$ does not satisfy Property (3). In this case
Property (4) shows that the least common multiple of the lengths of
all $u_1$-cycles and the least common multiple of the lengths of all
$u_j$-cycles in $\Gamma^{\ast}_{k+1}$ are different for some
$j=2,...,m$. So if $\varphi_1$ is the homomorphism of the group $G$
which corresponds to the action of $G$ on the set of vertices of the
graph $\Gamma_{k+1}^{\ast}$ then $\varphi_1(G)$ is a finite group
and $\varphi_1(u_1),...,\varphi_1(u_m)$ are nonidentity $p$-elements
and $\{u_1,..., u_m\}=\alpha\cup\beta$, where $\alpha$ and $\beta$
are nonempty sets with empty intersection, moreover, for each
$u_i\in\alpha$ and $u_j\in\beta$, we have $|\varphi_1(u_i)|>
|\varphi_1(u_j)|$.

If $m=2$ the proof of this statement shows that lemma 3 is true for
the case $m=2$. So we may apply the induction on $m$.

The set $\{u_1, ..., u_m\}$ can be presented as the union of two
disjoint subsets $\{z_1, ..., z_s\}$ and $\{y_1, ..., y_r\}$ and,
for each $i$ and $j$ such that $1\leqslant i\leqslant s$ and
$1\leqslant j\leqslant r$, we have $|\varphi_1(z_i)|>
|\varphi_1(y_j)|$. Besides, $\varphi_1(u_1), ..., \varphi_1(u_m)$
are nonidentity $p$-elements.

Put $\pi'=\pi\cup p$. By the induction hypothesis, there exists a
homomorphism $\varphi_2$ of group $G$ onto a finite group such that
the images of $y_1, ..., y_r$ have pairwise different orders and
$p\nmid\ |\varphi_2(y_i)|$ for each $i$ and each $p\in \pi'$

Suppose that $\rho$ is the set of all prime divisors of the numbers
$|\varphi_2(u_1)|$,\dots, $|\varphi_2(u_m)|$, $\pi''=\pi'\cup\rho$.
By the induction hypothesis, there exists a homomorphism $\varphi_3$
satisfying the requirements for elements $\{z_1, ..., z_s\}$ and the
set of primes $\pi''$. The homomorphism $\varphi:
G\rightarrow\varphi_1(G)\times\varphi_2(G)\times\varphi_3(G)$
defined by the formula $\varphi: g\mapsto(\varphi_1(g),
\varphi_2(g), \varphi_3(g))$ is as required.

\textbf{ Lemma 4.}\textsl{ Suppose that $G=A\ast B$, where the
groups $A$ and $B$ are finite, $w$ is a nonidentity cyclically
reduced element of the Cartesian subgroup $C$, and the integers
$k_1, ..., k_l$ have pairwise different absolute values. Then there
exists a homomorphism $\varphi$ of $G$ onto a finite group such that
the orders of the images of the elements $w^{k_1}, ..., w^{k_l}$ are
pairwise different.}

\textsl{ Proof. }

For each prime divisor $p_i$ of $k_1...k_l$, take a positive integer
$n_i$ such that $p_i^{n_i}$ does not divide $k_1...k_l$. According
to Lemma 1, there exists a homomorphism $\varphi_{p_i}$ of $G$ onto
a finite group such that $\varphi_{p_i}(w)$ is a $p_i$-element and
$|\varphi_{p_i}(w)|>p_i^{n_i}$. Let $p_1, ..., p_s$ be the set of
prime divisors of the numbers $k_1, ..., k_l$. Then the homomorphism
$\varphi:
G\rightarrow\varphi_{p_1}(G)\times...\times\varphi_{p_s}(G)$ defined
by the formula $\varphi : g\mapsto\prod_{i=1}^s\varphi_{p_i}(g)$ is
as required.

\section{Proof of Theorems.}

For each homomorphisms $\varphi$ and $\psi$ of groups $A$ and $B$
respectively we define the homomorphism $\varphi\ast\psi$ of $A\ast
B$ by the following way. If $s=s_1...s_k\in A\ast B$ is the normal
form of $s$, then $\varphi\ast\psi(s)=\phi_1(s_1)...\phi_k(s_k)$,
where $\phi_i=\varphi$ in case $s_i\in A$ and $\phi_i=\psi$ if
$s_i\in B$.

\textsl{ Proof of Theorems 1 and 2.}

Suppose that $G=A\ast B$ where $A$ and $B$ are either infinite
cyclic groups or periodic $n$-element order separable for each $n$
groups such that $|a|\neq|b|$ for each $a\in A\setminus\{1\}$ and
$b\in B\setminus\{1\}$.

Consider the elements $u_1, ..., u_n$ of $A\ast B$ such that $u_i$
and $u_j^{\pm1}$ are not conjugate for $i\ne j$

Without loss of generality, we assume that $u_1, ..., u_n$ are
cyclically reduced. Consider the following decomposition of the set
$\{u_1, ..., u_n\}=\alpha\cup\beta\cup\gamma$, where the sets
$\alpha$, $\beta$, and $\gamma$ are disjoint. The set $\alpha$
consists of all elements of the group $A$, $\beta$ consists of all
elements of the group $B$, $\gamma$ consists of all elements not
belonging to $A\cup B$.

We shall prove that there exist homomorphisms $\varphi$ and $\psi$
of $A$ and $B$ respectively onto finite groups such that the orders
of the images of elements of the set $\alpha\cup\beta$ are pairwise
different, each element of the set $\alpha\cup\beta\setminus\{1\}$
does not belong to the kernel of $\varphi$ or $\psi$, and the images
of elements of $\gamma$ do not belong to groups conjugate to
$\varphi(A)$ or $\psi(B)$ after the homomorphism $\varphi\ast\psi$
of the group $G$, and the elements $\{\varphi\ast\psi(u_1), ...,
\varphi\ast\psi(u_n)\}$ are such that $\varphi\ast\psi(u_i)$ and
$\varphi\ast\psi(u_j^{\pm1})$ are not conjugate for $i\ne j$.

Consider the set
$$
\Lambda= \alpha\cup\Omega\cup(\Omega\cdot\Omega^{-1}),
$$
where $\Omega$ is the set of all elements of the group $A$ occurring
in the normal forms of the elements of $\gamma$ ( if
$\Omega=\{\omega_j| j\in J\}$ then
$\Omega\cdot\Omega^{-1}=\{\omega_j\omega_k^{-1}| j, k\in J\}$).
Next, let us choose a maximal subset $\Lambda'$ of $\Lambda$ such
that for each different $x$ and $y$ from $\Lambda'$ $x$ and
$y^{\pm1}$ are not conjugate. Note that $1\in\Lambda'$.

Consider a homomorphism $\varphi$ of the group $A$ onto a finite
group such that the orders of the images of elements from $\Lambda'$
are pairwise different.

Similarly, we consider a homomorphism $\psi$ of $B$ for the
corresponding elements of this group. Thus, we may assume that $A$
and $B$ are finite ( we may consider that for each $a\in\alpha,
b\in\beta$, the  elements $\varphi\ast\psi(a)$ and
$\varphi\ast\psi(b)$ have different orders, because of the
properties of the free factors). If the set $\gamma$ is empty use
the residual finiteness of the group $A\ast B$. Suppose that
$\gamma$ is not empty. We may also assume that the set $\gamma$ lies
in the Cartesian subgroup. Indeed, $|G : C|< \infty$. Hence, there
exists a positive integer $l$ such that $u^l\in C$ for each
$u\in\gamma$. If there exists a homomorphism such that the orders of
the images of $u^l$ are pairwise different for each $u$ from
$\gamma$, then the same is true for the images of the elements $u$.
Furthermore, if $x,y \notin A\cup B$, $x$ and $y$ are cyclically
reduced and nonconjugate, then $x^n$ and $y^n$ are nonconjugate for
each $n$ (Combinatorial group theory,\ \textsl{Springer-Verlag},
theorem 1.4). Let us decompose the image of the set $\gamma$ into a
union of disjoint subsets $\alpha_1, ..., \alpha_s$, where
$\alpha_i$ is a maximal subset in $\gamma$ such that each element of
$\alpha_i$ belongs the cyclic subgroup conjugate to the subgroup
generated by some element $w_i$, and $w_i$ is the same for each
element of $\alpha_i$. We may assume that the set $\alpha_i$
consists of elements $w_i^{k_{ij}}$, and $w_i$ is not a proper
power. Let $d_i$ be the greatest common divisor of numbers $k_{ij}$
where $i$ is fixed. Then $w_i^{d_i}\in C$. There exists the number
$m_i$ such that $m_i|d_i, w_i^{m_i}\in C$ and if $m_i>1$ then for
each divisor $m_i'$ of $m_i$ which is less than $m_i$
$w_i^{m_i'}\notin C$. Thus we may consider that the set $\alpha_i$
consists of elements $w_i^{k_{ij}}$, $w_i\in C$ and if
$w_i=w_i'^{m_i}$, where $w_i'$ is not a proper power then
$w_i'^{m_i'}\notin C$ for each $m_i'$ such that $m_i'<m_i,
m_i'|m_i$. By Lemma 4, for each $\alpha_i$, there exists a
homomorphism $\varphi_i$ such that the orders of elements from
$\varphi_i(\alpha_i)$ are pairwise different.

The orders of the images of the elements from $\alpha\cup\beta$
after $\varphi_i$ are pairwise different because $\varphi_i\mid_A,
\varphi_i\mid_B$ are injective.

The set $\alpha_i$ consists of elements $w_i^{k_{ij}}$. Let $\pi$ be
the set consisting of prime divisors of the orders of elements
$\varphi_i(\alpha_i)$ and of the images of $\alpha$ and $\beta$
under the homomorphism
$G\rightarrow\varphi_1(G)\times...\times\varphi_s(G)$ defined by
$g\mapsto(\varphi_1(g), ..., \varphi_s(g))$. By Lemma 3 there exists
a homomorphism $\psi$ such that the orders of elements $\psi(w_1),
..., \psi(w_s)$ are pairwise different and $p\nmid\ |\psi(w_i)|$ for
each $i$ and for each $p\in\pi$ and $\psi\mid_A, \psi\mid_B$ are
injective.

The homomorphism $\chi: G \rightarrow\varphi_1(G)\times
...\times\varphi_s(G)\times\psi(G)$ defined by $\chi :
g\mapsto(\varphi_1(g), ..., \varphi_s(g), \psi(g))$ is as required.

\textsl{ Proof of Theorem 3.}

It is sufficient to show, how to reduce the general case to the
case, when the free factors are finite. The proof of this particular
case is similar to the proof of Theorems 1 and 2, because after we
reduced the general case to this particular case we did not use any
properties of finite free factors.

Suppose the elements $u, v, w$ are cyclically reduced and, if $x$
and $y$ are two different elements from $\{u,v,w\}$, then $x$ and
$y^{\pm1}$ are not conjugate.

If none of these elements lies in $A$, or none of them lies in $B$,
then the proof is the same as the proof of Theorems 1 and 2. Thus,
we have to consider two cases: $u\in A\setminus e, v\in B, w\notin
A\cup B$ or $u, v\in A, w\in B\setminus e$. In the first case, there
exists a homomorphism $\varphi_1$ of $A$ onto a finite group such
that $\varphi_1(u)\neq e$ and the image of $w$ does not belong to a
group conjugate to the image of some free factor after the
corresponding homomorphism of the group $G$, besides, if the order
of $u$ is infinite and the order of $v$ is finite, we may assume
that $\varphi_1(u^p)\neq e$ for each divisor $p$ of $|v|$. There
exists a homomorphism $\varphi_2$ of $B$ onto a finite group such
that $\varphi_2(v^{|\varphi_1(u)|})\neq e$ if $v\neq e$. In
addition, the image of $\varphi_1\ast$id$(w)$ does not belong to a
group conjugate to the image of a free factor after the
corresponding homomorphism of $\varphi_1(A)\ast B$. If $v=e$, then
it is sufficient to satisfy the first condition. In this case, we
use the residual finiteness of $A$ and $B$.

In the second case, if $u\neq e, v\neq e$ there exists a
homomorphism $\varphi_1$ of $A$ onto a finite group such that
$|\varphi_1(u)|\neq |\varphi_1(v)|$ and $\varphi_1(u)\neq e,
\varphi_1(v)\neq e$. The homomorphism $\varphi_1$ with the above
properties exists because the group $A$ is 3-element order
separable, and the natural homomorphism $\varphi_2:\varphi_1(A)\ast
B\rightarrow\varphi_1(A)$ is as required. Suppose that $u=e$. If the
orders of $v$ and $w$ are finite then $|v|$ and $|w|$ are coprime
and we use the residual finiteness of $G$. Consider the case $u=e$
and the order of $v$ is infinite. Since the group $A$ is residually
finite there exists a homomorphism $\varphi_1$ of $A$ onto a finite
group such that $\varphi_1(v)\neq e$. And if the order of $w$ is
finite we may consider that $\varphi_1(v^p)\neq e$ for each divisor
$p$ of $|w|$. The group $B$ is residually finite so there exists a
homomorphism $\varphi_2$ of $B$ onto a finite group such that
$\varphi_2(w^{|\varphi_1(v)|})\neq e$, and we use the residual
finiteness of the group $\varphi_1(A)\ast\varphi_2(B)$.

\begin{center}
\large{Acknowledgements}
\end{center}

The author thanks A. A. Klyachko for valuable comments and setting
the problem.

\begin{center}
\Large{ References }
\end{center}

1. Endimioni, G. Elements with the same normal closure in a
metabelian group. \textsl{The Quarterly Journal of Mathematics}
\textbf{58(1)}: 23-29, 2006.

2. Kargapolov, M. I., Merzlyakov, Yu. I. (1977). Foundations of group
theory. \textsl{Nauka}.

3. Klyachko,\ A. A. Equations over groups, quasivarieties, and a
residual property of a free group. \textsl{Journal of group theory}
\textbf{2}: 319-327, 1999.

4. Lyndon, R. C., Schupp, P. E. (1977). Combinatorial group
theory. \textsl{Springer-Verlag}.

5. Yedynak,\ V. V. Separability with respect to order. \textsl{Vestnik
Mosk. Univ. Ser. I Mat. Mekh.} \textbf{3}: 56-58, 2006.

Yedynak V. V., a graduate student of Moscow State Univ.


Yedynak V. V., Trudovaya street 14-b, 32, Ivanteevka, Moscow Region,
141281, Russia.

\end{document}